\documentclass[10pt]{article}

\usepackage{a4wide}
\usepackage{amssymb}
\usepackage{amsfonts}
\usepackage{amsmath}
\input xy
\xyoption{arrow} \xyoption{matrix}

\date{}

\newtheorem{proposition}{Proposition}[section]
\newtheorem{theorem}[proposition]{Theorem}
\newtheorem{lemma}[proposition]{Lemma}

\def\Kdim{{\rm K.dim }\,}

\def\der{\partial }

\def\nFM0{{\nu }_{F,M_0}}
\def\nFN0{{\nu }_{F,N_0}}
\def\nGN0{{\nu }_{G,N_0}}

\def\N0{ {\bf N}_0 }

\def\t{\otimes}
\def\g{\gamma}

\def\ra{\rightarrow}

\def\Xpm{X^{\pm }}

\def\s{\sigma}

\def\l1{{\lambda}_1}

\def\a{\alpha}
\def\a0{ {\alpha }_0}
\def\a1{ {\alpha }_1}

\def\l{\lambda}
\def\o{\omega}

\def\nFGM0{{\nu }_{F,G,M_0}}

\def\nFN0{{\nu}_{F,N_0}}

\def\sm{{\sigma}^m}

\def\sm1{{\sigma}^{-1}}

\def\smtp1{{\sigma}^{-t+1}}

\def\o{\omega }
\def\S1{S^{-1}}

\def\Xpm1{X^{\pm 1}_1}

\def\sPM1{{\sigma }^{\pm 1}}
\def\sMP1{{\sigma }^{\mp 1 }}

\def\d{\delta}

\def\di{{\rm d.ind}}

\def\L{\Lambda}

\def\G{\Gamma}

\def\CA{{\cal A}}

\def\CD{{\cal D}}

\def\Ytm1{Y^{t-1}}
\def\Yim1{Y^{i-1}}

\def\ass{{\rm ass}}

\def\Aut{{\rm Aut}}
\def\bK{\overline{K}}

\def\Der{{\rm Der }}

\def\gcd{ {\rm gcd } }

\def\SL2Z{ {\rm SL}_2({\bf Z}) }

\def\Gp1{ G^{1 , 1 } }
\def\P11{ P^{-1 , 1 } }
\def\Pp1{ P^{1 , 1 } }

\def\CE{{\cal E}}

\def\nCLsr{{}^\nu\kern-2pt {\cal L}^{\sigma , \rho  }}
\def\nP{{}^\nu \kern-2pt P}
\def\nL{{}^\nu\kern-2pt L}
\def\nLL{{}^\nu\kern-2pt \Lambda}
\def\nPsr{{}^\nu\kern-2pt P^{\sigma , \rho  }}
\def\nLsr{{}^\nu\kern-2pt L^{\sigma , \rho  }}
\def\nuCL{{}^\nu\kern-2pt  {\cal L}}
\def\nCLsr{{}^\nu\kern-2pt {\cal L}^{\sigma , \rho  }}
\def\nCL1m{{}^\nu\kern-2pt {\cal L}^{-1 , 1  }}
\def\x1nu{x^\frac{1}{\nu}}
\def\xm1nu{x^{-\frac{1}{\nu}}}

\def\ra{\rightarrow }

\def\CB{{\cal B}}

\def\CC{ {\cal C}}
\def\CE{ {\cal E} }

\def\nAM0{{\nu }_{{\cal A},M_0}}
\def\nAN0{{\nu }_{{\cal A},N_0}}

\def\Kdim{ {\rm Kdim } }

\def\Der{ {\rm Der }}

\def\gp{\mathfrak{p}}

\def\Spec{{\rm Spec}}

\def\di!{\frac{\der^i}{i!}}
\def\dik!{\frac{\der^k_i}{k!}}

\def\CC{{\cal C}}

\def\gldim{{\rm gldim}}

\def\N{\mathbb{N}}

\def\0{\overline{0}}
\def\1{\overline{1}}

\def\Ln1{\L_{n,\overline{1}}}

\def\a1{a_{\overline{1}}}

\def\St{{\rm St}}
\def\S{\Sigma}

\def\vn1{\overrightarrow{n-1}}

\def\mL{\mathbb{L}}

\def\mS{\mathbb{S}}
\def\mJ{\mathbb{J}}
\def\mI{\mathbb{I}}

\def\mU{\mathbb{U}}

\def\K1{{\rm K}_1}

\def\hmI1{\widehat{\mI_1}}
\def\tmI1{\widetilde{\mI_1}}
\def\tmJ1{\widetilde{\mJ_1}}
\def\hB1{\widehat{B_1}}
\def\hCB1{\widehat{\CB_1}}

\def\mL{\mathbb{L}}

\def\CB{{\cal B}}
\def\CC{{\cal C}}

\def\hB{\hat{B}}

\def\Alg{{\rm Alg}}

\begin{document}

\author{V. V.   Bavula 
} 
 
\title{Isomorphism problems and groups  of automorphisms for  Ore extensions $K[x][y; \d ]$ (zero characteristic)}

\maketitle 
\begin{abstract} 
Let $\L (f) = K[x][y; f\frac{d}{dx} ]$ be an Ore extension of a  polynomial  algebra $K[x]$ over a field $K$ of characteristic zero where $f\in K[x]$.  For a given polynomial $f$, the automorphism group of the algebra $\L (f) $ is explicitly described. The polynomial case $\L (0) = K[x,y]$ and the case of the Weyl algebra $A_1=  K[x][y; \frac{d}{dx} ]$ were done  done by Jung (1942) and van der Kulk (1953),  and  Dixmier  (1968), respectively. In 1997, Alev and Dumas proved that  the algebras $\L (f)$ and $\L (g)$ are isomorphic iff $g(x) = \l f(\alpha x+\beta )$ for  some $\l, \alpha \in K\backslash \{ 0\}$ and $\beta\in K$. In 2015, Benkart, Lopes and Ondrus gave a complete description of  the set of automorphism groups  of algebras $\L(f)$. In this paper we complete the picture, i.e. {\em given}  the polynomial $f$ we have the explicit description of the automorphism group of $\L (f)$.

The key concepts in finding the automorphism groups are the eigenform, the eigenroot  and the eigengroup of a polynomial (introduced in the paper; they are of independent interest). 
 \\



\end{abstract}


\section{Introduction}\label{INTR}

In this paper, $K$ is a field of characteristic zero and $\bK$ is  its algebraic closure, $K^\times := K\backslash \{ 0\}$,  $K[x]$ is a polynomial algebra in the variable  $x$ over $K$, $\Der_K(K[x])=K[x]\frac{d}{dx}$ is the set of all $K$-derivations of the algebra $K[x]$, 
$$\L :=\L (f):=K[x][y; \d := f\frac{d}{dx} ]=K\langle x, y \, |  \, yx-xy=f\rangle =\bigoplus_{i\geq 0} K[x]y^i$$
is an Ore extension of the algebra $K[x]$ where $f=f(x)\in K[x]$. Given  an algebra $D$ and its derivation $\d$, the {\em Ore extension} of $D$, denoted $D[y; \d]$,  is an algebra generated by the algebra $D$ and $y$ subject to the defining relations $yd-dy=\d (d)$ for all $d\in D$.  The algebra $\L$ is a Noetherian domain of Gelfand-Kirillov dimension 2.\\

By dividing the element $y$ by the leading coefficient of the polynomial $f$ we can assume that the polynomial $f$ is {\em monic}, i.e. its leading coefficient is 1 provided $f\neq 0$. Then the  algebras $\{ \L (f)\, | \, f\in K[x]\}$ as a class can be divided into four subclasses: $f=0$, $f=1$, the polynomial $f$ has only a {\em single} root in $\bK$ and the polynomial $f$ has at least two {\em distinct} roots in $\bK$.\\

 If $f=0$ then the algebra $\L (0)=K[x,y]$ is a polynomial algebra in two variables and its group of automorphisms is well-known \cite{Jung-1942,van der Kulk-1953}:
  The group $\Aut_K (K[x,y])$ is generated by the automorphisms:
 \begin{eqnarray*}
 t_l &:&  x\mapsto \l x, \;\;  \;\; \;\; \;\;  \;\; \; y\mapsto y, \\
\Phi_{n,\l} &:& x\mapsto x+\l y^n, \;\; y\mapsto y,  \\
 \Phi_{n,\l}' &:& x\mapsto x, \;\;  \;\; \;\; \;\;  \;\; \;\; \;\; y\mapsto y+\l x^n,
\end{eqnarray*}
 where $n\geq 0$ and $\l \in K$.

  If $f=1$ 
   then 
 the algebra $\L (1)$ is the (first) {\em Weyl algebra} $$A_1=K\langle 
x, \der\, | \, \der x-x\der =1\rangle\simeq K[x][y; \frac{d}{dx} ].$$
 In 1968, Dixmier \cite{Dix} gave an explicit generators for  the automorphism group $\Aut_K (A_1)$: The group $\Aut_K (A_1)$ is generated by the automorphisms:
 \begin{eqnarray*}
\Phi_{n,\l} &:& x\mapsto x+\l y^n, \;\; y\mapsto y,  \\
 \Phi_{n,\l}' &:& x\mapsto x, \;\, \; \;\; \;\;  \;\; \;\; \;\; y\mapsto y+\l x^n,
\end{eqnarray*}
 where $n\geq 0 $ and $\l \in K$. \\

Extending results of Dixmier \cite{Dix}  on the automorphisms of the Weyl algebra $A_1$,  
Bavula and Jordan \cite{Bav-Jor} considered isomorphisms and automorphisms of generalized
Weyl algebras over polynomial algebras of characteristic 0. Explicit generators of the automorphism group were given (they are more involved comparing to the case of the Weyl algebra to present them here). Alev and Dumas \cite{Alev-Dumas-1997}
initiated the study of automorphisms of Ore extensions $\L (f)$ in characteristic zero case. Their results were extended also to prime characteristic by Benkart, Lopes and Ondrus \cite{Benkart-Lopes-Ondrus-2015}. The  algebra $\L (x^2)$ (the Jordan plane)  was studied by Shirikov \cite{Shirikov-2015},  Cibils, Lauve, and Witherspoon \cite{CLW}, and  Iyudu \cite{Iyudu-2014}. Gadis \cite{Gaddis-2015} studied isomorphism problems for  algebras on two generators that satisfy a single quadratic relation.

 In 1981, a classification of simple $A_1$-modules was obtained by Block (over the field  of complex numbers) in  
 \cite{Bl} (see also \cite{Bav-SimMod-1992, Bav 2, Bav 3} for an alternative approach  via generalized Weyl algebras in a more general situation and over an arbitrary field).

In \cite{Bav 5, Bav-SimModOreExt} a classification of simple $\L (f)$-modules is given. \\

Benkart, Lopes and Ondrus  \cite{Benkart-Lopes-Ondrus-2015} (Theorems 8.3 and 8.6) gave a 
description of the {\em set} of  automorphisms groups  of algebras $\L (f)$ over arbitrary fields and if the  automorphism  group of $\L (f)$ is {\em given} they presented  information  on the type of  the polynomial $f$, \cite[Corollary 8.7]{Benkart-Lopes-Ondrus-2015} (in general, if one fixes the type of the polynomial then the automorphism group is {\em larger} than the one which  is naively expected). In this paper, we proceed in the opposite direction: if the polynomial $f$ is {\em given} then the automorphism group $\Aut_K\, \L (f)$   is explicitly described. \\


 
 
 {\bf Isomorphism problems  for  the algebras $\L (f)$.} Theorems \ref{Alev-Dum-P3.6}  is an isomorphism criterion for the algebras $\L$. It also describes the automorphism group of each algebra $\L (f)$.  
 
 \begin{theorem}\label{Alev-Dum-P3.6}
\cite[Proposition 3.6]{Alev-Dumas-1997} Let $f, g\in K[x]$ be   polynomials. Then $\L (f) \simeq \L(g)$ iff $g(x)= \l f(\alpha x+\beta)$ for some elements $\l , \alpha \in K^\times$ and $\beta \in K$.
\end{theorem}

{\bf The group of automorphisms of  the algebra $\L (f)$.}  The automorphism group $\Aut_K(\L (f))$ of the algebra $\L (f)$ contains an obvious subgroup 
\begin{equation}\label{grmS}
\mS :=\mS (K):= \{ s_p\, | \, p\in K[x]\}\simeq (K[x],+), \;  s_p\mapsto p\;\; {\rm where}\;\;  s_p(x) =x\; {\rm and } \;  s_p(y)=y+p.
\end{equation}
For each polynomial $f\in K[x]$,  we introduce a group  $G_f$ (the  eigengroup group of $f$, see below).  Lemma \ref{a5Mar20} describes an explicit monomorphism $G_f\ra \Aut_K(\L (f))$. We identify the group $G_f(K)$ with its image in $ \Aut_K(\L (f))$. Theorem \ref{5Mar20} follows from Theorem \ref{Alev-Dum-P3.6}. 

\begin{theorem}\label{5Mar20}
 Suppose that $f\in K[x]$ is a monic nonscalar polynomial. Then $$\Aut_K(\L (f))=\mS (K) \rtimes G_f(K).$$ 
\end{theorem}
So, the theorem above states that the group $\Aut_K(\L (f))$ is a semidirect product of its  two subgroups ($\mS (K)$ is a normal subgroup of $\Aut_K(\L (f))$).  \\

{\bf The eigengroup $G_f(K)$ of a polynomial $f\in K[x]$.} Given a group $G$,  a $G$-module $V$ over a field $K$ and a non-empty subset $U$ of $V$. The {\em eigengroup} of the set $U$ in $G$, denoted by $G_U(K)$, is the set of all elements of the group $G$ such that the elements of the set $U$ are common eigenvectors of  them with eigenvalues in the  field $K$. Clearly, the eigengroup is a subgroup of $G$ and $$G_U=\bigcap_{u\in U} G_u$$ where $G_u:=G_{\{ u\} }=\{g\in G\, | \, gu=\l(g)u$ for some $\l (g)\in K\}$. If $K$ is a subfield of a field $K'$ then $G_U(K)\subseteq G_U(K')$ where $U$ is a subset of the $G$-module $K'\t_KV$ over the field $K'$.

We are interested in the case when $G=\Aut_K(K[x])=\{ \s_{\l , \mu} \, | \, \l \in K^\times, \mu \in K\}$ and $V=K[x]$ where $\s_{\l , \mu } (x) = \l x+\mu$. Clearly, $G_\mu =\Aut_K(K[x])$ for all $\mu \in K$ and $G_f=G_{\nu f}$  for all elements $f\in K[x]$ and $\nu\in K^\times$. So, in order to describe the eigengroups for a nonscalar polynomial we can assume that it is a  monic  polynomial.  \\

{\em Definition.} Let $f(x) = x^d+a_{d-1}x^{d-1}+\cdots+a_1x+a_0\in \bK [x]$ be a monic polynomial of degree $d\geq 1$ where $a_i\in \bK$ are the coefficients of the polynomial $f(x)$. Then the natural number
$$ \gcd (f(x)):=\gcd \{ i\geq 1\, | \, a_i\neq 0\}$$
is called the {\em exponent} of $f(x)$.\\

Clearly, the exponent of $f(x)$ is the largest natural number $m\geq 0$ such that $f(x)=g(x^m)$ for some polynomial $g(x)\in K[x]$. \\

{\em Definition.} Let $f(x) = x^d+a_{d-1}x^{d-1}+\cdots+a_1x+a_0\in \bK [x]$ be a monic polynomial of degree $d\geq 1$ where $a_i\in \bK$ are the coefficients of the polynomial $f(x)$. The polynomial $f(x)$ admits a {\em unique} presentation, which we call the {\em igenform} or {\em eigenpresentation} of $f(x)$,
\begin{equation}\label{eigenf}
f(x) = (x-\nu )^s g\Big((x-\nu)^n\Big)
\end{equation}
where
\begin{itemize}
\item  $\nu = \nu_f :=-\frac{a_{d-1}}{d}\in K$,
\item  $s=s_f\geq 0$ is the multiplicity of the factor $x-\nu$ (i.e. if $\nu $ is a root of $f(x)$ then $s$ is its multiplicity and $s=0$ otherwise), 
\item  $n=n_f:=\gcd \bigg(\frac{f(x+\nu)}{x^s}\bigg)$, and
\item   $g(x)=g_f(x) \in \bK [x]$ is a monic polynomial such that $d=s+n\deg (g(x))$.
\end{itemize}
 The scalar $\nu\in \bK$ and the natural number $s$ are  called the {\em eigenroot}  of $f(x)$ and  its {\em multiplicity}, respectively. The natural number $n\geq 1$ and the monic polynomial $g(x)$ are called the {\em eigenorder}  and the {\em eigenfactor} of $f(x)$. In general, the eigenroot $\nu$ may not be a root of the polynomial $f(x)$. Clearly, $g(0)\neq 0$. \\

The uniqueness of the eigenform of the polynomial $f(x)$ is obvious from the definition.  Notice that the multiplicity $s_f$ is the unique natural  number $s\geq 0$ such that $$f(\nu ) = f'(\nu) = \cdots = f^{(s-1)}(\nu ) =0\;\; {\rm  and }\;\; f^{(s)}(\nu ) \neq 0$$ where $f^{(i)}(x):=\frac{d^if(x)}{dx}$ is the $i$'th derivative of the polynomial $f(x)$.\\ 

If the polynomial $f(x)=(x-\nu')^d$ has only single root $\nu'\in K$ then 
$f(x)=(x-\nu')^d$ is the eigenform of $f(x)$, i.e. $\nu_f = \nu'$, $s=d$, $n=0$    and $g(x) =1$. \\

For each nonscalar monic polynomial $f(x)$, Theorem \ref{EG5Mar20} and Theorem \ref{AEG5Mar20} are explicit descriptions of the eigengroup $G_f(K)$ in the case when the field $K$ is algebraically closed and in general case, respectively.

\begin{theorem}\label{EG5Mar20}
Suppose that $K = \bK$ is an algebraically closed field,  $f\in K[x]$ is a monic polynomial of degree $d\geq 1$  and $f(x) = (x-\nu )^s g\Big((x-\nu)^n\Big)$ is its eigenform. 
\begin{enumerate}
\item If  the polynomial $f(x)=(x-\nu)^d$ has only single root $\nu\in K$ then $G_f=\{ \s_{\l , (1-\l ) \nu}\, | \, \l \in K^\times \} \simeq K^\times$,  $\s_{\l , (1-\l ) \nu}\mapsto \l $ where $\s_{\l , (1-\l ) \nu} (x)=\l x+(1-\l ) \nu$. 
\item Suppose that  the polynomial $f(x)$ has at least two distinct roots. Then 
\begin{enumerate}
\item $G_f\neq \{ e\}$ iff  $n>1$.  
\item If $n>1$ then $G_f=\langle \s_{\l_n,(1-\l_n ) \nu} \rangle = \{   \s_{\l_n,(1-\l_n ) \nu}^j\, | \, 0\leq j \leq n-1\}$ is the cyclic group of order $n$ where $ \l_n := e^\frac{2\pi i}{n}$ (is the  primitive $n$'th root of unity  and $i=\sqrt{-1}$) and $\s_{\l_n , (1-\l_n ) \nu} (x)=\l_n x+(1-\l_n ) \nu$. 
\end{enumerate}
\end{enumerate}
\end{theorem}

\begin{theorem}\label{AEG5Mar20}
Suppose that $K$ is not necessarily algebraically closed field,  $f= x^d+a_{d-1}x^{d-1}+\cdots + a_0\in K[x]$ is a monic polynomial of degree $d\geq 1$  and $f(x) = (x-\nu )^s g\Big((x-\nu)^n\Big)$ is its eigenform (as an element of  $\bK [x]$). Then $\nu \in K$ and  $g(x) \in K[x]$.
\begin{enumerate}
\item  $G_f(K)=G_f(\bK ) \cap \Aut_K(K[x])=\{ \s \in G_f(\bK ) \, | \, \s (x)\in K[x]\}$. 
\item If  the polynomial $f(x)=(x-\nu)^d$ has only single root $\nu$  then $G_f=\{ \s_{\l , (1-\l ) \nu}\, | \, \l \in K^\times \} \simeq K^\times$,  $\s_{\l , (1-\l ) \nu}\mapsto \l $ where $\s_{\l , (1-\l ) \nu} (x)=\l x+(1-\l ) \nu$. 
\item Suppose that  the polynomial $f(x)$ has at least two distinct roots in $\bK$. Then 
\begin{enumerate}
\item $G_f\neq \{ e\}$ iff  $n>1$ and $\l_m:= e^\frac{2\pi i}{m}\in K$ for some natural number $m\geq 2$ such that $m|n$.   
\item If $n>1$ and $n':= \max \{ m \, | \,  m|n, \,\l_m\in K\}\geq 2$ 
then $$G_f=\langle \s_{\l_{n'},(1-\l_{n'} ) \nu} \rangle = \{   \s_{\l_{n'},(1-\l_{n'} ) \nu}^j\, | \, 0\leq j \leq n'-1\}$$ is the cyclic group of order $n'$ where $ \l_{n'} := e^\frac{2\pi i}{n'}$  and $\s_{\l_{n'} , (1-\l_{n'} ) \nu} (x)=\l_{n'} x+(1-\l_{n'} ) \nu$. 
\end{enumerate}
\end{enumerate}
\end{theorem}

{\bf The inclusion $G_f\subseteq \Aut_K(\L (f))$.} 

\begin{lemma}\label{a5Mar20}
Suppose that  $f= x^d+a_{d-1}x^{d-1}+\cdots + a_0\in K[x]$ is a monic polynomial of degree $d\geq 1$  and $f(x) = (x-\nu )^s g\Big((x-\nu)^n\Big)$ is its eigenform (as an element of $\bK [x]$).
\begin{enumerate}
\item If  the polynomial $f(x)=(x-\nu)^d$ has only single root $\nu$  then the map $$G_f(K)=\{ \s_{\l , (1-\l ) \nu}\, | \, \l \in K^\times \}\ra \Aut_K(\L (f)), \;\; \s_{\l , (1-\l ) \nu}\mapsto \s_{\l , (1-\l ) \nu}$$ is a group monomorphism where $\s_{\l , (1-\l ) \nu} (x)=\l x+(1-\l ) \nu$ and $\s_{\l , (1-\l ) \nu} (y)=\l^{d-1} y$.
\item Suppose that  the polynomial $f(x)$ has at least two distinct roots in $\bK$ and $G_f(K)\neq \{ e\}$, i.e. $n>1$ and $n':= \max \{ m \, | \,  m|n, \, \l_m\in K\}\geq 2$, by Theorem \ref{AEG5Mar20}.(3). Then the map 
$$G_f(K)=\langle \s_{\l_{n'},(1-\l_{n'} ) \nu} \rangle\ra \Aut_K(\L (f)),\;\;  \s_{\l_{n'},(1-\l_{n'} ) \nu}\mapsto\s_{\l_{n'},(1-\l_{n'} ) \nu}$$ is a group monomorphism where $\s_{\l_{n'} , (1-\l_{n'} ) \nu} (x)=\l_{n'} x+(1-\l_{n'} ) \nu$ and $\s_{\l_{n'} , (1-\l_{n'} ) \nu} (y)=\l_{n'}^{d-1} y$.  
\end{enumerate}
\end{lemma}

{\it Proof}. The proof follows at once from the fact that the image of the automorphisms respect the defining relation $[y,x]=f$ of the algebra $\L (f)$. $\Box $\\

Proofs of 
Theorem \ref{EG5Mar20} and   Theorem \ref{AEG5Mar20} are given in Section \ref{AUTR}. 
 In \cite{Bav-AutOreCharp}, similar results are obtained in prime characteristic and they are much more involved comparing to the characteristic zero case.


\section{Isomorphism problems and groups  of automorphisms for  Ore extensions $K[x][y; \d ]$}\label{AUTR}

We keep the notation of Section \ref{INTR}.\\

{\bf Proof of Theorem \ref{EG5Mar20}.} 1. Let $\s = \s_{\l ,\mu}$ where $\s_{\l , \mu } (x) =\l x+\mu$. Then $\s \in G_f$ iff $\s (x-\nu ) = \l (x-\nu )$ iff $\nu = \frac{\mu }{1-\l}$ iff $\s =\s_{\l , (1-\l ) \nu}$. 

2. Suppose that $n>1$.  Notice that $\s_{\l_n,(1-\l_n ) \nu}(x-\nu ) = \l_n (x-\nu)$. Hence, 
$\s_{\l_n,(1-\l_n ) \nu}(f)= \l_n^s f$, and so $\s_{\l_n,(1-\l_n ) \nu}\in G_f$.  For all $j\geq 1$, 
$$\s_{\l_n,(1-\l_n ) \nu}^j=\s_{\l_n^j,(1-\l_n^j )\nu }.$$
Therefore, $\langle \s_{\l_n,(1-\l_n ) \nu} \rangle$ is the cyclic group of order $n$ which is a subgroup of $G_f$. In particular, $G_f\neq \{ e\}$. 

Conversely, suppose that $G_f\neq \{ e\}$ and $e\neq \s_{\l , \mu} \in G_f$.
 Since $\s_{\l , \mu } (f) = \l^df$, the automorphism $\s_{\l , \mu }$ permutes the roots of the polynomial $f$ (i.e. the minimal primes of the ideal $(f)$ of the polynomial algebra $K[x]$). Since 
 $$\s_{\l , \mu}^j=\s_{\l^j, (1+\l +\cdots +\l^{j-1})\mu}\;\; {\rm  for \; all}\;\; j\geq 0$$ and every minimal prime over $(f)$ is equal to the ideal $(x-\xi)$ where $\xi$ is a root of the polynomial $f$, we must have $\s_{\l , \mu}^m=e$ for some natural number $m\geq 2$ (since the polynomial $f$ has at least two distinct roots). Hence, $\l^m=1$. 
Clearly, $m\leq d!$. Let $\mU :=\{ \g \in K\, | \, \g^t=1$ for some $t\geq 1\}$ be the group of all roots of unity in $K$. We have the group homomorphism
 $$ \phi : G_f \ra \mU, \;\; \s_{\l , \mu} \mapsto \l$$
 since $\s_{\l ,\mu} \s_{\l', \mu'}= \s_{\l \l', \l'\mu +\mu'}$ for all $\l , \l'\in K^\times$ and $\mu , \mu'\in K$. Since the order of elements in the image, say $\G_f$, of the homomorphism above are bounded by $d!$ the group 
$\G_f$ is equal to $\langle \l_{n'}\rangle$ where $\l_{n'}=e^\frac{2\pi i }{n'}$ for some natural number $n'\geq 2$. We will see that $n'=n=n_f$. 
 The homomorphism $\phi$ is a monomorphism since otherwise $\s _{1, \rho} \in G_f$ for some $\rho\in K^\times$ but this is impossible as the order of the element $\s _{1, \rho} $ is infinite since $$\s _{1, \rho}^j=\s _{1, j\rho}\; \; {\rm  for\;  all} \;\; j\geq 1$$ and the field $K$ has characteristic zero. Therefore, the group $G_f\simeq \G_f$ is a cyclic group of order $n'$ which is generated by the element $\phi^{-1}(\l_{n'}) = \s_{\l_{n'}, (1-\l_{n'})\nu'}$ for some element $\nu'\in K$.  Let $\s':=\s_{\l_{n'}, (1-\l_{n'})\nu'}$.
 
 Since $\s'(x-\nu')=\l_{n'}(x-\nu')$, 
 $$K[x]=\bigoplus_{j=0}^{n'-1}(x-\nu')^jK[x]^{\s'}\;\; {\rm where} \;\; K[x]^{\s'}=K[(x-\nu')^{n'}]$$
 is the fixed  ring of the automorphism $\s'$, i.e. $K[x]^{\s'}:=\{ g\in K[x]\, | \, \s'(g)=g\}$. Therefore, 
$$f(x) = (x-\nu' )^{s'} g'\Big((x-\nu')^{n'}\Big)$$
 for some natural number $s'\geq 0$  and a polynomial $g'\in K[x]$ such that $g'(0)\neq 0$. Clearly, the sum of the roots of the  polynomial $f$ is equal to $-a_{d-1}=d\nu'$ where $f= x^d+a_{d-1}x^{d-1}+\cdots + a_0$ and $a_j$ are the coefficients of the polynomial $f$. Then 
 $$\nu'=-\frac{a_{d-1}}{d}=\nu_f=\nu .$$
 Hence,  $s'=s$ and $n'\leq n$, by the very definition of the element $n=n_f$. On the other hand, $\s_{\l_n, (1-\l_n)\nu}\in G_f=\langle \s_{\l_{n'}, (1-\l_{n'})\nu'} \rangle$, and so $n|n'$, i.e. $n=n'$. The proof of statement 2 is complete.    $\Box $\\

{\bf Proof of Theorem \ref{AEG5Mar20}.} Clearly, $\nu =-\frac{a_{d-1}}{d}\in K$. Then $x-\nu \in K[x]$, and so  and $g(x) \in K[x]$, by the very definition of the polynomial $g(x)$.

1. Statement 1 is obvious.
 
2. Statement 2 follows from statement 1, Theorem \ref{EG5Mar20}.(1) and the fact that $\nu \in K$.

3. Statement 3 follows from statement 1 and Theorem \ref{EG5Mar20}.(2). $\Box$ \\

 {\bf The prime spectrum of the algebra $\L $.} An element $a$ of an algebra $A$ is called a {\em normal element} of $A$ if $Aa=aA$. An element $a$ of an algebra $A$ is called a {\em regular element} if it is not a zero divisor. The set of all regular elements of the algebra $A$ is denoted by $\CC_A$. An ideal $\gp$ of a ring $R$ is called a {\em completely prime ideal} if the factor ring $R/\gp$ is a domain. A completely prime ideal is a prime ideal. The sets of prime and completely prime ideals of the ring $R$ are denoted by $\Spec (R)$ and $\Spec_c(R)$, respectively. 

\begin{theorem}\label{KGLD}
(\cite[Theorem 1.1]{Bav-SimModOreExt}) Let $K$ be a field of characteristic zero, $\L =K[x][y; \d := f\frac{d}{dx} ]$ where $f\in K[x]\backslash K$. Let $f=p_1^{n_1}\cdots p_s^{n_s}$ be a unique (up to permutation) product of irreducible polynomials of $K[x]$. Then 
\begin{enumerate}
\item The Krull dimension of $\L$ is  $\Kdim (\L )=2$.
\item The global  dimension of $\L$ is  $\gldim(\L )=2$.
\item The elements $p_1, \ldots , p_s$ are regular normal elements of the algebra $\L$ (i.e. $p_i$ is a non-zero-divisor of $\L$ and $p_i\L = \L p_i$). 
\item $\Spec (\L )=\Spec_c (\L )=\{ 0 , \L p_i,  (p_i, q_i)\, | \, i=1, \ldots , s;\; q_i\in  {\rm Irr}_m(F_i[y])\}$  where $F_i:=K[x]/(p_i)$ is a field  and ${\rm Irr}_m(F_i[y])$ is the set of monic irreducible polynomials of  the polynomial algebra $F_i[y]$ over the field $F_i$ in the variable $y$. If, in addition, the field $K$ is an algebraically closed and $\l_1, \ldots , \l_s$ are the roots of the polynomial $f$ then $\Spec (\L )=\{ 0 , \L (x-\l_i),  (x-\l_i, y-\mu)\, | \, i=1, \ldots , s;\; \mu\in K\}$. 
\end{enumerate}
\end{theorem}

Given an action of a  group  $G$ on a set $S$ and an element $s\in S$. The set $\St_G(s):=\{ g\in G\, | \, gx=x\}$ is called the {\em stabilizer} of the element $s$ in $G$.  The stabilizer is a subgroup of $G$, and the map $G/\St_G(s)\ra Gs$, $g\St_G(s)\mapsto gx$ is a bijection where $G/\St_G(s) := \{ g\St_G(s)\, | \, g\in G\}$ is the set of right $\St_G(s)$-cosets and $Gs$ is the $G$-orbit of  the element $s$.

Given a group $G$, a normal subgroup $N$ and a subgroup $H$. The  group $G$ is called the {\em semidirect product} of $N$ and $H$, written $G=N\rtimes H$,  if $G=NH:=\{ nh\, | \, n\in N, h\in H\}$ and $N\cap H=\{ e\}$ where $e$ is the identity of the group  $G$.  \\

{\bf Proof of Theorem \ref{5Mar20}.} Let $\s$ be an automorphism of the $K$-algebra $\L = \L (f)$. It can be uniquely extended to a $\bK$-automorphism, say $\s$,  of the algebra $\bK\t_K\L$. Let $\l_1, \ldots , \l_s$ be the roots of the polynomial $f$ in $\bK$. By Theorem \ref{KGLD}.(4), the automorphism $\s$ permutes the set $\{ (x-\l_1), \ldots , (x-\l_s)\}$ of height 1 prime ideals.  Since the elements $x-\l_1, \ldots , x-\l_s$ are regular normal elements of the domain $\bK \t_K\L$ and the set $\bK^\times$ is the group of units of the algebra $\bK\t_K \L$, we must have that 
$$\s (x) = \l x+\mu $$ for some elements $\l\in \bK^\times$ and $\mu \in \bK$.   Since $K[x]=\L \cap \bK [x]$, we must have that $\s (x)\in \s (\L ) \cap \s (\bK [x] )= \L \cap \bK [x] =K[x]$, and so  $\l \in K^\times$ and $\mu \in K$.  So, the automorphism $\s$ respects the polynomial algebra $K[x]$. In particular it respects the Ore set $S=K[x]\backslash \{ 0\}$ of the algebra $\L$. The automorphism $\s$ can be uniquely extended to  an automorphism of the algebra $B_1=S^{-1}\L =K(x)[\der ; \frac{d}{dx}]$. Then $\s (\der ) = \l^{-1}\der +q$ for some element $q\in K[x]$ (since $\s (K(x))=K(x)$). In particular, 
$$\s (y) = \s (f\der )= \s (f) (\l^{-1} \der +q)=\l^{-1} \frac{\s (f)}{f} y +p\;\; {\rm where}\;\; p:=\s (f) q\in K[x] $$
and $\s (f) = \g f$ for some element $\g \in K^\times$. Clearly, $\g = \l^d$ where $d=\deg (f)$ is the degree of the polynomial $f$ (since $\s (x) = \l x+\mu$). So,
$$\s (x)=\l x+\mu \;\; {\rm and}\;\; \s (y)=\l^{d-1}y+p, $$
i.e. $\s \in \mS (K)\rtimes G_f(K)$, as required. $\Box$. \\

{\bf Isomorphism problems and groups  of automorphisms for subalgebras of an algebra that respect a left denominator set.} 
For an algebra $A$, we denote by $\CC_A$ its set of {\em regular} elements, i.e. non-zero divisors. A non-empty subset $S$ of the algebra $A$ is called a {\em multiplicative set} (or {\em multiplicatively closed set}) if $SS\subset S$, $1\in S$ and $0\not\in S$. A multiplicative  set $S$ of $A$ is called a {\em left Ore set} of $A$ if the {\em left  Ore condition} holds: For each pair of elements $s\in $ and $ a\in A$, there are elements $s'\in S$ and $a'\in A$ such that $s'a=a's$. If $S$ is a left Ore set of the algebra $A$ then the set $\ass (S):=\{ a\in A\, | \, sa=0$ for some $s\in S\}$ is an ideal of $A$. A left Ore set $S$ of $A$ is called a {\em left denominator set} of $A$ if $as=0$ for some elements $s\in S$ and $a\in A$ implies $s'a=0$  for some element $s'\in S$. If $S$ is a left denominator set of $A$ then the algebra $S^{-1}A:=\{ s^{-1}a\, \ \, s\in S, a\in A\}$ is called the  {\em localization} of $A$ at $S$. Clearly, is $S$ is a left Ore set that consists of regular elements then $S$ is a left denominator set with $\ass (S)=0$ and $A\subseteq S^{-1}A$. 

\begin{theorem}\label{2Mar20}
 Let $A$ be an algebra, $S$ be a left Ore  set of $A$  that consists of regular elements of $A$, $B:=S^{-1}A$, $R$ be a subalgebra of $A$ that contains $S$ as a left Ore set such that $S^{-1}R=B$ and $\CA (A,S) $ be the set of all such subalgebras $R$,  and  $\Aut_K(R, S):=\{ \s \in \Aut_K(R)\, | \, \s (S)=S\}$. Then 
\begin{enumerate}
\item $\Aut_K(R, S)=\{ \tau \in \Aut_K(B)\, | \, \tau (S)=S,  \tau (R)=R\}$.
\item Let $\CB (A,S):=\{ R\in \CA (A, S)\, | \,  \s (S)=S$  for all $\s \in \Aut_K(R) \}$. Then for every $R\in \CB (A, S)$, $\Aut_K(R)=\Aut_K(R, S)=\{ \tau \in \Aut_K(B)\, | \, \tau (S)=S,  \tau (R)=R\}$.
\item The algebras $R,R'\in \CB (A, S)$ are isomorphic iff there there is an automorphism $\s\in \Aut_K(B, S)$ such that   $\s (R) = R'$. 
The map 
$$\Aut_K(B; S, R) /\Aut_K(R, S)\ra \Alg (A, R), \;\; \s \Aut_K(R, S)\mapsto \s (R)$$
is a bijection where $\Aut_K(B; S, R):=\{ \s \in \Aut_K(B)\, | \, \s(S) = S, \; \s (R)\subseteq A\}$, $\Aut_K(B; S, R) /\Aut_K(R, S):= \{ \s \Aut_K(R, S)\, | \, \s \in \Aut_K(B; S, R)$ and $\Alg (A, R)$ is the set of all subalgebras $R'\in \CB (A, S)$ that are isomorphic to the algebra $R$. 

\item Let $\CC (A, S):=\{ R\in \CB (A, S)\, | \, \Aut_K(B; S, R)=\Aut_K(R, S)\}$. The the algebras in the set $\CC (A, S)$ are all non-isomorphic. 
\item Let $\CE (A, S):=\{ R\in\CC (A, S)\, | \, \Aut_K(R) =  \Aut_K(R, S)\}$. The the algebras in the set $\CE (A, S)$ are all non-isomorphic and  $\Aut_K(R, S)=\{ \tau \in \Aut_K(B)\, | \, \tau (S)=S,  \tau (R)=R\}$ for all $R\in \CE (A, S)$. 
\end{enumerate}
\end{theorem}

{\it Proof}. 1. By the very definition, the set  $\mL := \{ \tau \in \Aut_K(B)\, | \, \tau (S)=S,  \tau (R)=R\}$ is a subgroup of the group $\Aut_K(B)$ and the restriction map $\mL\ra \Aut_K(R, S)$, $\s \mapsto \s |_R$ is a group monomomorphism since $B=S^{-1}R$. We identify the group $\mL$ with its image in $\Aut_K(R, S)$, i.e. $\Aut_K(R, S)\supseteq \mL$.

It remains to prove that the reverse inclusion holds. Given an element $\s \in \Aut_K(R, S)$. We have to show that $\s \in \mL$. Since $\s (S)=S$, the automorphism $\s$ of the algebra $A$ can be {\em uniquely} extended to an automorphism, say $\s$, of the localization $B=S^{-1}R$:   For all elements $s\in S$ and $a\in A$, $\s (s^{-1} a) = \s (s)^{-1}\s (a)$. Since $\s (R)=R$, we have that $\s \in \mL$. 

2. Statement 2 follows from statement 1. 

3. Suppose that $\s : R\ra R'$ be an isomorphism of algebras  $R,R'\in \CB (A, S)$. Since $$R, R'\subseteq B=S^{-1}R= S^{-1}R'=(\s (S))^{-1}\s (R)= \s (S^{-1}R)=\s (B),$$ 
where that an automorphism $\s :B\ra B$ is the unique extension of the isomorphism $\s : R\ra R'$ (for all elements $s\in S$ and $r\in R$, $\s(s^{-1}r) = \s (s)^{-1}\s (r)$) and the diagram below is commutative: 

\begin{displaymath}
    \xymatrix{
        B \ar[r]^\s & B  \\
        R \ar[r]_{\s}\ar[u]       & R'\ar[u] }
\end{displaymath}
Clearly, $\s \in \Aut_K(B; S, R)$ and $\Aut_K(B; S, R)\Aut_K(R, S)=\Aut_K(B; S, R)$, and statement 3 follows. 

4. Statement 4 follows from statement 3.  

5. Statement 5 follows from statements 1 and  4. $\Box $\\

{\bf The algebra $B_1$ and its automorphism group.} Let $K(x)$ be the field of rational functions in the variable $x$. Then the Ore extension  $B_1:= K(x)[\der ; \frac{d}{dx}]$ is the localization $B_1=(K[x]\backslash \{ 0\})^{-1}A_1$ of the Weyl algebra $A_1$ at the Ore set $K[x]\backslash \{ 0\}$. The algebra $\L =\L (f)$ can be identified with a subalgebra of  the Weyl algebra $A_1$  by the monomorphism:
\begin{equation}\label{(2.1)}
\L \ra A_1,\;\;x\mapsto x,\;\; y\mapsto f\der.
\end{equation}
So, $\L = K\langle  x, y=f\der\rangle\subset A_1$. 
The element $f$ is a regular normal element of $\L$  (i.e. $\L f=f\L )$ since 
$$fy=yf-f'f=(y-f')f \;\;{\rm where}\;\;f'=\frac{df}{dx}.$$
 It determines the $K$-automorphism $\o_f$ of the algebra $\L $:
$$fu=\o_f(u)f,\;\; u\in \L ,$$
$$\o_f:x\mapsto x,\;\;y\mapsto y-f'.$$
We denote by $\L_f$  and  $A_{1,f}$ the localizations of the algebras  $\L $ and  $A_1$ at the powers of the element $f$, i.e.
$$\L_f=S^{-1}_f\L \;\;{\rm and  }\;\;A_{1,f}=S^{-1}_fA_1\;\;{\rm where}\;\;
S_f=\{ f^i\}_{i\geq 0}.$$
By (\ref{(2.1)}), 
\begin{equation}\label{(2.3)}
\L \subset A_1\subset \L_f= A_{1, f}=K[x, f^{-1}][\der ;\frac{d}{dx}]\subset B_1.
\end{equation}

The Weyl algebra $A_1=\CD (K[x])$ is the algebra of differential operators $\CD (K[x])$ on the polynomial algebra. The field of fractions $K(x)=(K[x]\backslash \{ 0\} )^{-1}K[x]$ of the polynomial algebra $J[x]$ at the Ore set $K[x]\backslash \{ 0\}$. Hence, the algebra 
$$B_1=(K[x]\backslash \{ 0\} )^{-1}A_1=(K[x]\backslash \{ 0\} )^{-1}\CD (K[x])\simeq \CD ((K[x]\backslash \{ 0\} )^{-1}K[x])=\CD (K(x))$$
is the ring of differential operators on the algebra $K(x)$. For any commutative algebra $A$, every automorphism $\s\in \Aut_K(A)$ can be uniquely extended to an automorphism $\s$ of the algebra $\CD (A)$ of differential operators on $A$ by the rule (the chain rule, the change of variables rule): For all $u\in \CD (A)$, $\s (u):= \s u\s^{-1}$. Therefore, $\Aut_K(A)\subseteq \Aut_K(\CD (A)$. In particular, 
 $$ \Aut_K(K[x])\subseteq \Aut_K(A_1)\;\; {\rm and}\;\; \Aut_K(K(x))\subseteq \Aut_K(B_1) $$
where $ \Aut_K(K[x])=\{ \s_{\l , \mu}\, | \, \l \in K^\times , \mu \in K\}$, $\s_{\l ,\mu}(x)=\l x+\mu$ and  $$ \Aut_K(K(x))=\{ \s_M\, | \, M\in {\rm PGL}_2(K)\}\simeq {\rm PGL}_2(K), \;\; \s_M\mapsto M\;\; {\rm  where}\;\;  
\s_M(x)=\frac{ax+b}{cx+d}, $$ 
$M=\begin{pmatrix}
a & b \\  c &d
 \end{pmatrix}\in  {\rm PGL}_2(K):={\rm GL}_2 (K)/K^\times E$ and $E=\begin{pmatrix}
1 & 0 \\ 0  &1
 \end{pmatrix}$.
The unique  extension of the automorphism  $\s_M\in \Aut_K(K(x))$ to an automorphism of the algebra $B_1$ is given by the chain rule: 
$$\s_M\frac{d}{dx} \s_M^{-1}=\bigg(\frac{dx'}{dx}\bigg)^{-1}\frac{d}{dx} \;\; {\rm where}\;\; x':=\s_M(x).$$

The automorphism group  $\Aut_K(B_1)$ acts in the obvious way on the algebra $B_1$. Let $\mS_1:= \St_{\Aut_K(B_1)}(x):=\{ \s\in \Aut_K(B_1)\, | \, \s (x)=x\}$, the stabilizer of the element $x\in B_1$ in $\Aut_K(B_1)$.

\begin{lemma}\label{a2Mar20}
\begin{enumerate}
\item $\Aut_K(B_1) = \mS_1\rtimes \Aut_K(K(x))=\{ \s_{M, q}\, | \, M\in {\rm PGL}_2(K), q\in K(x)\}$ where $ \s_{M, q}(x)=\s_M( x)$ ($\s_M$ is defined above) and  $ \s_{M, q}(\der )=\l^{-1}\der+q$;  $\mS_1=\{ s_q :=\s_{E,q}\, | \, q\in K(x)\}\simeq (K(x), +)$, $s_q\mapsto q$ where $s_q (x)=x$ and $s_q(\der )=\der +q$ ($E$ is the identity $2\times 2$ matrix).
\item  $\Aut_K(B_1, K[x]):=\{ \s \in \Aut_K(B_1)\, | \, \s (K[x])=K[x]\} = \mS_1\rtimes \Aut_K(K[x])=\{ \s_{l,\mu, q}\, | \, \l\in K^\times, \mu\in K, q\in K(x)\}$ where $ \s_{l,\mu, q}(x)=\l x+\mu$ and  $ \s_{\l,\mu, q}(\der )=\l^{-1}\der+q$.
\end{enumerate}
\end{lemma}

{\it Proof}. 1. Since  $\mS_1:= \St_{\Aut_K(B_1)}(x)$, we must have $\mS_1\cap \Aut_K(B_1)=\{ e\}$ and $\s \mS_1\s^{-1}\subseteq \mS_1$ for all automorphisms $\s\in \Aut_K(B_1)$. Hence,  $\Aut_K(B_1) \supseteq  \mS_1\rtimes \Aut_K(K(x))$.

To prove that the reverse inclusion holds we have to show that every  element $\s \in \Aut_K(B_1)$ belongs to the group $\mS_1\rtimes \Aut_K(K(x))$. The group of units  $K(x)^\times := K(x)\backslash \{ 0\}$ of the algebra $B_1$ is a $\s$-invariant set, i.e. $\s (K(x)^\times)=K(x)^\times$. Hence so is  the field $K(x)$. Let $\tau$ be the restriction of the automorphism $\s$ to the field $K(x)$. Then $\s_1:=\tau^{-1}\s\in \mS_1$, and so $\s = \tau \s_1\in \mS_1\rtimes \Aut_K(K(x))$, as required.

2. Statement 2 follows from statement 1. $\Box $\\

We give a different proof of Theorem \ref{Alev-Dum-P3.6} below.\\

{\bf Proof of Theorem \ref{Alev-Dum-P3.6}.} Let $\s : \L (f) \ra \L (g)$ be an isomorphism of the $K$-algebras. It can be uniquely extended to a $\bK$-isomorphism  $\s :\bK\t_K \L (f) \ra \bK\t_K\L (g)$. Let $\l_1, \ldots , \l_s$ (resp., $\l_1', \ldots , \l_t'$)  be the roots of the polynomial $f$ (resp., $g$)  in $\bK$. By Theorem \ref{KGLD}.(4), the automorphism $\s$ maps bijectively  the set $\{ (x-\l_1), \ldots , (x-\l_s)\}$ of height 1 prime ideals of the algebra $\bK\t_K \L (f)$ to the set $\{ (x-\l_1'), \ldots , (x-\l_t')\}$ of height 1 prime ideals of the algebra $\bK\t_K \L (g)$.  Therefore, $s=t$.  Since the elements $x-\l_1', \ldots , x-\l_t'$ are regular normal elements of the domain $\bK \t_K\L (g)$ and the set $\bK^\times$ is the group of units of the algebra $\L (g) $, we must have that 
$$\s (x) = \l x+\mu $$ for some elements $\l\in \bK^\times$ and $\mu \in \bK$.   Since $K[x]=\L (g) \cap \bK [x]$, we must have that $\s (x)\in \s (\L (f)) \cap \s (\bK [x] )= \L (g) \cap \bK [x] =K[x]$, and so  $\l \in K^\times$ and $\mu \in K$.  So, the isomorphism $\s$ respects the polynomial algebra $K[x]$ of the algebras $\l (f)$ and $\L (g)$. In particular it respects the Ore sets $S=K[x]\backslash \{ 0\}$ of the algebras $\L (f)$ and $\L (g)$, i.e. $\s (S)=S$. The isomorphism $\s$ can be uniquely extended to  an automorphism of the algebra $\s : B_1=S^{-1}\L (f) \ra B_1=S^{-1}\L (g)$. Then $\s (\der ) = \l^{-1}\der +q$ for some element $q\in K[x]$. In particular, 
$$\s (y) = \s (f\der )= \s (f) (\l^{-1} \der +q)=\l^{-1} \frac{\s (f)}{g} y +p\;\; {\rm where}\;\; p:=\s (f) q\in K[x] $$
and $\s (f) = \g g$ for some element $0\neq \g \in K[x]$. Applying the same argument for the isomorphism $\s^{-1}: \L (g) \ra \L (f)$, we have that $\s^{-1} (g) = \g_1 f$ for some element $0\neq \g_1 \in K[x]$.
Therefore, $f=\s^{-1}\s (f)=\s^{-1}(\g g)=\s^{-1}(\g ) \g_1f$, and so $\g , \g_1\in K^\times$, $\g_1=\g^{-1}$. 
Clearly, $\L (f) = \L (g)$ and $\g = \l^d$ where $d=\deg (f)$ is the degree of the polynomial $f$ (since $\s (x) = \l x+\mu$). So,
$$\s (x)=\l x+\mu \;\; {\rm and}\;\; \s (y)=\l^{d-1}y+p, $$
and the theorem follows. $\Box$. 

$${\bf Acnowledgements} $$

The author would like to thank the Royal Society for support.

Department of Pure Mathematics

University of Sheffield

Hicks Building

Sheffield S3 7RH

UK

email: v.bavula@sheffield.ac.uk


\begin{thebibliography}{6} 

\bibitem{Alev-Dumas-1997}  J. Alev and F. Dumas, Invariants du corps de Weyl sous l’action de groupes finis.
{\em Comm. Algebra} {\bf 25} (1997) 1655–1672.


\bibitem{Bav-SimMod-1992} V. V. Bavula, Simple $D[X,Y;\sigma,a]$-modules. {\it  Ukrainian Math. J.} {\bf 44} (1992) no. 12, 1500--1511.

\bibitem{Bav 2} V. V. Bavula, Generalized Weyl algebras and their representations,  {\em Algebra i Analiz} {\bf  4} (1992), no. 1, 75-97; English transl. in  {\em St.Petersburg Math. J.} {\bf 4} (1993), no. 1, 71--92.

\bibitem{Bav 3} V. V. Bavula, Generalized Weyl algebras, kernel and tensor-simple algebras, their simple modules, Representations of algebras. Sixth International Conference, August 19--22, 1992. CMS Conference proceedings ( V.Dlab and H.Lenzing Eds.), v. {\bf 14} (1993), 83--106. 


\bibitem{Bav 5} V. V. Bavula,  The simple modules of the Ore extensions with coefficients from   Dedekind ring. {\it  Comm.  Algebra}  {\bf 27} (1999) no. 6, 2665--2699.

\bibitem{Bav-Jor} V. V. Bavula and D. A. Jordan, Isomorphism problems and groups of
   automorphisms for generalized Weyl algebras. {\em Trans. Amer. Math. Soc.}  {\bf 353} (2001) no. 2,
   769--794.

\bibitem{Bav-SimModOreExt}  V. V. Bavula,  Classification of simple modules of the Ore extension $K[X][Y; f\frac{d}{dX}]$, {\em Mathematics in Computer Sciences} {\bf 14} (2020) 317--325. 

\bibitem{Bav-AutOreCharp} V. V. Bavula, Isomorphism problems and groups  of automorphisms for  Ore extensions $K[x][y; f\frac{d}{dx} ]$ (prime characteristic), submitted. 


\bibitem{Benkart-Lopes-Ondrus-2015}
G. Benkart, S. A. Lopes and  M.  Ondrus, A parametric family of subalgebras of the Weyl algebra I. Structure and automorphisms. {\em Trans. Amer. Math. Soc.} {\bf 367} (2015), no. 3, 1993–2021.

\bibitem{Bl} R. E. Block, The irreducible representations of the Lie algebra $sl(2)$ and of the Weyl algebra,  {\em Adv. Math.}  {\bf 39} (1981) 69--110. 

\bibitem{Gaddis-2015} J. Gaddis,  Two-generated algebras and standard-form congruence. {\em  Comm. Algebra} {\bf 43} (2015), no. 4, 1668--1686.

\bibitem{CLW} C. Cibils, A. Lauve, and S. Witherspoon, Hopf quivers and Nichols algebras in positive characteristic, {\em Proc. Amer. Math. Soc.} {\bf  137} (2009), no. 12, 4029--4041, 

\bibitem{Dix}
J. Dixmier,  Sur les alg\`{e}bres de Weyl,
{\it Bull. Soc. Math. France} {\bf 96} (1968), 209--242.

\bibitem{Iyudu-2014} N. K. Iyudu, 
Representation spaces of the Jordan plane. {\em  
Comm. Algebra} {\bf  42} (2014), no. 8, 3507--3540. 

\bibitem{Jung-1942}  H. W. E. Jung, \"{U}ber ganze birationale Transformationen der Ebene, {\em J. Reine Angew. Math.} {\bf 184} (1942), 161--174.

\bibitem{Shirikov-2015} E. N. Shirikov, 
Two-generated graded algebras. 
{\em Algebra Discrete Math.} (2005), no. 3, 60--84.

\bibitem{van der Kulk-1953} W. van der Kulk, On polynomial rings in two variables, {\em Nieuw Arch. Wiskunde} {\bf 1} (1953), 33--41.


























\end{thebibliography}
\end{document}